\documentclass[12pt]{article}%21.04.2016
\usepackage{amsmath}
\usepackage{amssymb}
\usepackage{dsfont}
\usepackage{cite}
\newsymbol \blackbox 1004
\newcommand{\eh}{\hfill}\newlength{\sperr}

\newenvironment{proof}{{\settowidth{\sperr}{\bf\rm
Proof}%
\par\addvspace{0.3cm}\noindent\parbox[t]{1.3\sperr}
{\bf\rm P\eh r\eh o\eh o\eh f\eh }%
}}{\nopagebreak\mbox{}
$\blackbox$\par\addvspace{0.3cm}}

\def\nn{\nonumber}
\def\a{\alpha}

\def\Lam{\Lambda}
\def\s{\sigma}
\def\la{\lambda}
\def\om{\omega}

\def\wh{\widehat}
\def\wt{\widetilde}
\def\ov{\overline}

\def\p{\partial}

\def\BC{{\mathbb C}}

\def\BN{{\mathbb N}}

\def\cli{{\mathcal I}}

\def\clm{{\mathcal M}}

\def\clu{{\mathcal U}}

\def\col{\mathrm{col}}
\def\diag{\mathrm{diag}}
\newcommand{\E}{\mathrm{e}}
\newcommand{\I}{\mathrm{i}}

\newtheorem{Pa}{Paper}[section]
\newtheorem{Tm}[Pa]{{\bf Theorem}}

\newtheorem{Cy}[Pa]{{\bf Corollary}}
\newtheorem{Rk}[Pa]{{\bf Remark}}

\newtheorem{Pn}[Pa]{{\bf Proposition}}

\title{Hamiltonian systems and \\ Sturm-Liouville equations: \\  Darboux transformation and applications}

\author{Alexander Sakhnovich}

\date{}
\begin{document}
\maketitle

\begin{abstract} We introduce GBDT version of Darboux transformation for
symplectic and Hamiltonian systems as well as for Shin-Zettl
systems and Sturm-Liouville equations. These are the first results on
Darboux transformation for general-type Hamiltonian and for Shin-Zettl
systems. The obtained results are applied to the corresponding transformations of the
Weyl-Titchmarsh functions and to the construction of explicit solutions 
of dynamical symplectic systems, of  two-way diffusion equations
and of indefinite Sturm-Liouville equations.
The energy of the explicit solutions of dynamical systems is expressed
(in a quite simple form) in terms of the parameter matrices of GBDT.
\end{abstract}



\section{Introduction}\label{Intro}
\setcounter{equation}{0}
This paper is dedicated to the study of the important subclasses
of the first order differential systems with a spectral parameter $\la$.
Namely, we consider Hamiltonian systems
\begin{align}& \label{E0}
\frac{d}{dx}y(x,\la)= F(x,\la)
y(x,\la), \quad F(x,\la)=J\big(\la H_1(x)+H_0(x)\big),
\end{align}
where
\begin{align}& \label{E0!}
J^*=-J, \quad H_1(x)^*=H_1(x), \quad H_0(x)=H_0(x)^*, \quad H_1(x)\geq 0;
\end{align}
and so called Shin-Zettl systems
\begin{align}& \label{E1}
\frac{d}{dx}y(x,\la)= F(x,\la)
y(x,\la), \quad F(x,\la)=\begin{bmatrix} r_1(x) & p(x)^{-1} \\ q(x) - \la \om(x) & \, r_2(x)
\end{bmatrix}.
\end{align}
Here $J^*$  is the conjugate transpose of the matrix $J$.
We assume that the $m \times m$ ($m\in \BN$) matrix functions $H_1(x)$ and $H_0(x)$ in \eqref{E0} and the functions 
$p^{-1}$, $q$, $r_1$, $r_2$ and $\om$ in \eqref{E1} are locally summable on $[0,\ell)$ ($\ell \leq \infty$). 
The matrix function $F$ in \eqref{E1} is the $2 \times 2$ Shin-Zettl matrix
of general form (see, e.g., $\S$ 2 in  \cite{Ev1} or in \cite{Ev2}). We note that Shin-Zettl differential expressions
were introduced in \cite{Shin, Ze} and were actively studied in regularization and  spectral theories (see the books \cite{AmP, Z}, papers \cite{Ev1,Ev2},
recent  surveys  \cite{Mirz, ZeS} and various references therein). The Lagrange-symmetric case
\begin{align}& \label{E4}
\om= \ov \om, \quad p=\ov p, \quad q=\ov q, \quad r_1=-\ov{r_2}
\end{align}
and the Lagrange-J-symmetric case
\begin{align}& \label{E5}
 r_1=-r_2
\end{align}
are of special interest \cite{Ev2}.
Here $\ov{\mu}$  stands for the value which is complex conjugate to $\mu$.

The entries of the $2\times 1$ vector function $y$ in \eqref{E1} are
denoted by $y_1$ and $y_2$. 
When $r_1 \equiv r_2 \equiv 0$, we rewrite \eqref{E1} in the
form
\begin{align}& \label{E2}
y_1^{\prime}=p^{-1}y_2, \quad y_2^{\prime}=(q-\la \om)y_1 \quad \left(y_k^{\prime}=\frac{d}{dx}y_k\right),
\end{align}
which is equivalent to the {\it Sturm-Liouville equation}
\begin{align}& \label{E3}
-\big(p(x) u^{\prime}(x,\la)\big)^{\prime}+ q(x)u(x,\la)=\la \om(x) u(x,\la),
\end{align}
where $u=y_1$. If $\om=\ov \om$, $p =\ov p$ and $\om$ or $p$ change signs, one speaks about {\it  indefinite  Sturm-Liouville problem}.
Quasi-derivatives related to the quasi-derivatives  generated by Shin-Zettl systems are used  in the study of important modifications of
Schr\"odinger-type operators (see, e.g.,  \cite{EGNST, ScTr} and references therein) including Schr\"odinger-type operators with distributional  potentials \cite{EGNST}.

On the other hand, Lagrange-symmetric Shin-Zettl systems, where $\om~\geq~0$, form also a subclass of Hamiltonian systems.
See, for instance, \cite{HiSh} on the representation \eqref{E0}, \eqref{E0!} of Hamiltonian systems and the equivalence of the definite Sturm-Liouville equation
to a certain subclass of Hamiltonian systems. We note that the book \cite{Atk} by Atkinson, the papers
by Hinton and Shaw as well the Kac-Krein supplement \cite{KaKr} (to the translation of \cite{Atk}) presented seminal developments in the
theory of Hamiltonian systems and Sturm-Liouville equations. (For recent references on Hamiltonian systems see, e.g., \cite{JaZ, Mog, SaSaR, SeSi}.)
 In some works, conditions \eqref{W1} are added in the definition of Hamiltonian systems but these conditions
are absent in \cite{HiSh} and they are not essential for Darboux transformations, which we will construct here, as well.

In this paper we construct our GBDT version of the B\"acklund-Darboux transformation (see  the results
and references in \cite{SaA2, SaA6, SaSaR}) for the cases of Hamiltonian and Shin-Zettl systems in order to study
perturbations  of these systems and corresponding transformations of
the Weyl-Titchmarsh functions. We construct explicit solutions of the perturbed systems as well.
Several versions of B\"acklund-Darboux transformations (see, e.g., \cite{Ci, Gu, MS, SaSaR} and references
therein) are a well-known tool for the construction of explicit solutions of linear and integrable nonlinear
equations. 
GBDT as well as Crum-Krein and commutation methods
(which are related to B\"acklund-Darboux transformations) are also essential in the study of Weyl-Titchmarsh theory and important spectral
problems \cite{Cr, D, Ge, GeT, GKS6, KoSaTe, Kr, MST, SaA8}. 

As far as we know, neither B\"acklund-Darboux transformations
nor commutation methods were applied to general-type Hamiltonian systems \eqref{E0} and to Shin-Zettl systems \eqref{E1} before
(although commutation and B\"acklund-Darboux transformations for  such important particular cases as Schr\"odinger equations,
canonical systems and related Dirac equations are well-known).
 We mention an interesting paper \cite{BeH}
on Kummer-Liouville transformation for Shin-Zettl systems but that transformation is different and was applied with different
purposes.

Darboux transformation for symplectic and general-type Hamiltonian systems is introduced in Section \ref{sec1'}.
The  corresponding transformations of the Weyl-Titchmarsh functions are considered in Section \ref{secAp}.
GBDT for Shin-Zettl systems and Sturm-Liouville equations is introduced in Sections \ref{sec2}-\ref{StL}.
Explicit solutions of dynamical symplectic systems and of  two-way diffusion equations are constructed
in Section \ref{DS}. Finally, explicit solutions of indefinite Sturm-Liouville equations are considered in Section \ref{IStL}.

As usual, $\BN$ denotes the set of natural numbers, $\BC$ denotes the complex plane, $\BC_+$ is the open upper
half-plane $\{\la: \, \Im(\la)>0\}$ and $\BC_-$ is the open lower
half-plane $\{\la: \, \Im(\la)<0\}$.
The notation $I_n$ stands for the $n\times n$ identity matrix, $H^*$  is the conjugate transpose of the matrix $H$,
the inequality $H \geq 0$ means that $H=H^*$ and that all the eigenvalues of the matrix $H$ are nonnegative.

 \section{GBDT  for  Hamiltonian systems}\label{sec1'}
 \setcounter{equation}{0}
 \paragraph{1.}
 Our GBDT version of B\"acklund-Darboux transformation for system \eqref{E0} is a particular case of GBDT
 for systems with rational dependence on spectral parameter (see, e.g., \cite{SaA6} or \cite[Sect. 7.2]{SaSaR}).
 We start with introducing GBDT for  the  first order system of $m$  differential equations
 with  a linear dependence on the spectral parameter ($m \in \BN$):
 \begin{align}& \label{H1}
y^{\prime}(x,\la)= F(x,\la)
y(x,\la), \quad F(x,\la)=-\big(\la Q_1(x)+Q_0(x)\big).
\end{align}
For that purpose we fix some initial system \eqref{H1} (i.e., some $m \times m$  matrix functions $Q_1(x)$ and $Q_0(x)$,
which are  locally summable on $[0,\ell)$),
 an integer $n \in \BN$ and five parameter matrices, namely, $n \times n$ matrices $A_1$, $A_2$ and $S(0)$, and
 $n \times m $ matrices $\Pi_1(0)$ and  $\Pi_2(0)$ such that the matrix identity
 \begin{align}& \label{E6}
A_1S(0)-S(0)A_2=\Pi_1(0)\Pi_2(0)^*
\end{align}
 holds. Matrix functions $\Pi_1(x)$, $\Pi_2(x)$ and $S(x)$
are introduced by their initial values $\Pi_1(0)$, $\Pi_2(0)$, $S(0)$ and differential equations
 \begin{equation} \label{E7}
\Pi_1^{\prime}=A_1\Pi_1 Q_1+\Pi_1 Q_0, \quad (\Pi_2^*)^{\prime}=-Q_1\Pi_2^* A_2-Q_0\Pi_2^*, \quad S^{\prime}=\Pi_1 Q_1 \Pi_2^*.
\end{equation}
The identity 
\begin{align}& \label{E9}
A_1S(x)-S(x)A_2=\Pi_1(x)\Pi_2(x)^*,
\end{align}
for all $x\in [0,\ell)$, is a particular case of \cite[f-la (7.18)]{SaSaR} and easily follows from \eqref{E6} and \eqref{E7}.

{\it When we deal with $S(x)^{-1}$, our further statements are valid in the points of invertibility of $S(x)$.}
The questions of invertibility of $S(x)$ are discussed in our sections separately (see, e.g., Remarks \ref{RkS} and \ref{S>0}).

According to the subcase $r=1$, $\, l=0$ of \cite[Theor. 7.4]{SaSaR}, the so called Darboux matrix 
for system \eqref{H1} is given by the formula
\begin{align}& \label{E10}
w_A(x,\la)=I_m-\Pi_2(x)^*S(x)^{-1}(A_1-\la I_n)^{-1}\Pi_1(x).
\end{align}
More precisely, \cite[Theor. 7.4]{SaSaR} yields that $w_A$ satisfies the following
equation
\begin{align}& \label{H2}
\frac{d}{d x}w_A(x,\la)=\wt F(x,\la)w_A(x,\la)-w_A(x,\la) F(x,\la),
\end{align}
where
\begin{align}& \label{E11}
\wt F(x,\la):=-\big(\la Q_1(x)+ \wt Q_0(x)\big),  \\ 
& \label{E11!}
 \wt Q_0(x):=Q_0(x)- \big(Q_1(x)X(x)-X(x)Q_1(x)\big), \\
& \label{E12}
 X(x):=\Pi_2(x)^*S(x)^{-1}\Pi_1(x).
\end{align}
We note that (in view of \eqref{E9}) the matrix function $w_A(\la)$ of the form \eqref{E10}
is (for each $x$) the so called transfer matrix function in Lev Sakhnovich form
(see \cite{SaL1, SaL3, SaSaR} and references therein).

System $y^{\prime}=\wt Fy$ is called the transformed (GBDT-transformed) system (recall that \eqref{H1}
is the initial system).  An important step in the proof of \eqref{H2} is the proof of the equation
\begin{align}& \label{D0}
\big(\Pi_2^*S^{-1}\big)^{\prime}=-Q_1\Pi_2^*S^{-1}A_1-\wt Q_0\Pi_2^*S^{-1}.
\end{align}
See \cite[f-la (7.61)]{SaSaR} for the general formula, of which \eqref{D0} is a particular case.
We shall use \eqref{D0} in Section \ref{DSympl}.

Formula \eqref{H2} implies the following theorem.
\begin{Tm}\label{Tm1GBDT} Let $y(x,\la)$ satisfy  system \eqref{H1} and let
$w_A$ be given by \eqref{E10}, where the matrix functions $\Pi_1$, $\Pi_2$ and $S$ are determined
by \eqref{E7} and identity \eqref{E6} holds. Then the function 
\begin{align}& \label{H3}
\wt y(x,\la):=w_A(x,\la)y(x,\la)
\end{align}
satisfies, in the points of invertibility of $S(x)$, another $($transformed$)$ first order system
\begin{align}& \label{H4}
\frac{d}{d x}\wt y(x,\la)= \wt F(x,\la)
\wt y(x,\la),
\end{align}
where $\wt F(x,\la)$ is given by \eqref{E11}--\eqref{E12}.
\end{Tm}

\paragraph{2.} The most important subcase of the considered above GBDT-transforma-tions is the subcase 
of the initial  system \eqref{H1} such that
\begin{align}& \label{H5}
Q_1(x)=-JH_1(x), \quad Q_0(x)=-JH_0(x), \\
& \label{H6}
 J^*=-J, \quad H_1(x)^*=H_1(x), \quad H_0(x)=H_0(x)^*.
\end{align}
 In that subcase we deal with system \eqref{E0}, where all the conditions \eqref{E0!} on
 Hamiltonian system, excluding the nonnegativity condition $H_1(x) \geq 0$, hold.
 If $J^*=J^{-1}$ (e.g., $J$ has the form \eqref{W1})
 conditions \eqref{H6} mean that system \eqref{E0} is symplectic.
 Further in the paragraphs 2 and 3 we assume that the equalities \eqref{H5} and \eqref{H6} are valid. 
 
 We
 omit indices in $A_1$ and $\Pi_1$ and set
 \begin{align}& \label{H7}
A=A_1, \quad \Pi=\Pi_1; \quad  A_2=A^*, \quad S(0)=S(0)^*, \quad \Pi_2(0)= - \Pi(0)J.
\end{align}
Using \eqref{H5}--\eqref{H7} we rewrite the first and second equations in \eqref{E7}, correspondingly,
in the forms
 \begin{equation} \nn
(-\Pi J)^{\prime}=-A (-\Pi J) H_1J-(-\Pi J) H_0 J, \quad (\Pi_2)^{\prime}=-A \Pi_2 H_1J-\Pi_2 H_0 J.
\end{equation}
Thus, the equations on $-\Pi J$ and on $\Pi_2$ coincide, and, in view of $\Pi_2(0)=- \Pi(0)J$ we obtain
$\Pi_2(x)\equiv -\Pi(x)J$. In this way, equations \eqref{E7} are reduced to the equations
 \begin{equation} \label{H9}
\Pi^{\prime}=-A\Pi J H_1(x)- \Pi J H_0(x), \quad  S^{\prime}(x)=\Pi J H_1(x)J^* \Pi(x)^*.
\end{equation}
Since we assume in \eqref{H7} that $S(0)=S(0)^*$, the second equation in \eqref{H9} yields $S(x)=S(x)^*$.
Thus, we have
\begin{align}& \label{H9'}
\Pi_2(x)\equiv -\Pi(x)J, \quad S(x)=S(x)^*.
\end{align}
Now, the matrix identity \eqref{E9} and Darboux matrix \eqref{E10} are rewritten in the form
\begin{align}& \label{H10}
AS(x)-S(x)A^*=\Pi(x)J \Pi(x)^*,
\\  & \label{H11}
w_A(x,\la)=I_m-J\Pi(x)^*S(x)^{-1}(A-\la I_n)^{-1}\Pi(x).
\end{align}
Moreover, using \eqref{H5} and the equalities $\Pi_2(x)\equiv -\Pi(x)J$ and $J^*=-J$, we rewrite \eqref{E11}--\eqref{E12} in the form
\begin{align}& \label{H12}
\wt F(x,\la)=J\big(\la H_1(x)+\wt H_0(x)\big), \quad \wt H_0(x)=H_0(x)+Z(x), \\
& \label{H12!}
 Z(x):=\Pi(x)^*S(x)^{-1}\Pi(x)JH_1(x)+H_1(x)J^*\Pi(x)^*S(x)^{-1}\Pi(x).
\end{align}
Formulas \eqref{H6}, \eqref{H12} and \eqref{H12!} imply that $\wt H_0=\wt H_0^*$, that is, $\wt F$ has the same
form as $F$. Hence, the next proposition follows from Theorem \ref{Tm1GBDT}.
\begin{Pn}\label{Cy1GBDT} Let $y(x,\la)$ satisfy  system \eqref{E0} $($such that \eqref{H6} holds$)$, and  let
a triple $\{A, S(0)=S(0)^*, \Pi(0)\}$ of parameter matrices satisfying \eqref{H10} at $x=0$ be given.
Introduce $w_A(x,\la)$ by \eqref{H11}, where the matrix functions $\Pi(x)$  and $S(x)$ are determined
by \eqref{H9}. 

Then the function $\wt y(x,\la)=w_A(x,\la)y(x,\la)$
satisfies, in the points of invertibility of $S(x)$, another $($transformed$)$  system of the same form as \eqref{E0}, namely,
\begin{align}& \label{H15}
\frac{d}{d x}\wt y(x,\la)= \wt F(x,\la)
\wt y(x,\la),
\end{align}
where $\wt F(x,\la)$ is given by \eqref{H12}, \eqref{H12!} and the equality $\wt H_0=\wt H_0^*$ holds.

If $H_1(x)\geq 0$ and $S(0)>0$, the systems \eqref{E0} and \eqref{H15} are Hamiltonian.
\end{Pn}
\begin{Rk}\label{RkS} If system \eqref{E0}, \eqref{H6} is Hamiltonian $($i.e., $H_1(x)\geq 0)$ and, in addition, the inequality $S(0)>0$ holds,
formula \eqref{H9} shows that $S(x)>0$ for all $x \in [0,\ell)$. Therefore, $S(x)$ is invertible on $[0, \ell)$.
In particular, it follows that the system \eqref{H15} is, indeed, Hamiltonian.

\end{Rk}
\paragraph{3.} If in the system \eqref{E0} we have $J=-J^*=-J^{-1}$ and $H_0\equiv 0$, we come to the important class of canonical systems.
See GBDT for canonical system and its applications to Weyl-Titchmarsh theory in \cite{SaA5}.
 For the case of Hamiltonian systems with invertible $J$ we can (similar to the case of canonical  systems) consider
 transformation slightly different from \eqref{H12}, \eqref{H12!}. More precisely, we introduce
 matrix functions $\wh w(x)$ and $v(x,\la)$ by the formulas
\begin{align}& \label{H12'}
\wh w^{\prime}(x)=-\wh w(x) J Z(x), \quad \wh w(0)=I_m; \quad v(x,\la)=\wh w(x)w_A(x,\la).
\end{align}
It is easy to see that $\wh w(x) J \wh w(x)^*=J$, and so
\begin{align}& \label{H13}
\wh w(x)^{-1}= J \wh w(x)^*J^{-1}=J^* \wh w(x)^*(J^*)^{-1}.
\end{align}
In view of Proposition  \ref{Cy1GBDT} and relations \eqref{H12'} and \eqref{H13}, if $y(x,\la)$ satisfies \eqref{E0}, then
the matrix function $\wh y(x,\la)=v(x,\la)y(x,\la)$ satisfies the system
\begin{align}& \label{E0+}
\frac{d}{d x}\wh y(x,\la)= \wh F(x,\la)
\wh y(x,\la), \quad \wh F(x,\la)=J\big(\la \wh H_1(x)+\wh H_0(x)\big),
\end{align}
where
\begin{align}& \label{E16}
\wh H_1=J^{-1}\wh w JH_1\wh w^{-1}=J^{-1}\wh w JH_1J^* \wh w^*(J^*)^{-1}=\wh H_1^*,
\\ & \label{E17}
\wh H_0= J^{-1}\wh w J(\wt H_0 -Z)\wh w^{-1}=J^{-1}\wh w JH_0J^* \wh w^*(J^*)^{-1}=\wh H_0^*.
\end{align}
In the special case $H_0=\I c J^{-1}$ ($c=\ov c$), the formula \eqref{E17} is simplified and we obtain
$\wh H_0 \equiv \I c J^{-1}$.

%%%%%%%%%%%%%%%%%%%%%%%%%%%%%%%%%%%%%%%%%%%%%%%%%%%%%%%%%%%%%
 %%%%%%%%%%%%%%%%%%%%%%%%%%%%%%%
  \section{Darboux transformations \\ of Weyl-Titchmarsh functions}\label{secAp}
 \setcounter{equation}{0}
 In his important paper \cite{Kra}, Krall introduced Weyl-Titchmarsh (or simply Weyl) $M(\la)$-functions 
 of Hamiltonian systems in the classical terms of ``Weyl circle" inequalities. 
Here, Weyl circles of system \eqref{E0} on the intervals $[0,\ell']$ $(\ell'<\ell)$ and the values  $\la$
in the upper half-plane 
$\la \in \BC_+$ (i.e., $\Im(\la)>0$)
are considered. The Weyl circles in the lower half-plane $\BC_-$
are treated in a quite similar way and we omit that case. 
 
Krall required that $m$ is even and that $J$ in \eqref{E0}
 has a special form:
 \begin{align}& \label{W1}
m=2r \quad (r \in \BN), \quad J=\begin{bmatrix}0 & I_r \\ -I_r  & 0
\end{bmatrix}.
\end{align}
In fact, Hamiltonian system in \cite{Kra} is written in a slightly different from \eqref{E0} way and our $J^*$ stands for $J$ in an equivalent
to \eqref{E0} system in \cite{Kra}. 
Rewriting correspondingly the inequality for  the Weyl circle (of  matrices $M(\la)$ with
$\la \in \BC_+$)  from \cite[p. 670]{Kra},
we obtain
  \begin{align}& \label{W2}
\I\begin{bmatrix}I_r & M(\la)^*
\end{bmatrix}Y(\ell',\la)^*JY(\ell',\la)\begin{bmatrix}I_r \\ M(\la)
\end{bmatrix}\leq 0.
\end{align}
Here $Y(x,\la)$ is the fundamental $m\times m$ solution of the Hamiltonian system \eqref{E0} (such that \eqref{E0!} and \eqref{W1} are valid),
normalized by the initial condition
 \begin{align}& \label{W3}
Y(0,\la)=E \quad (EJ=JE, \quad E^*E=I_m).
\end{align}

According to Proposition \ref{Cy1GBDT}, the fundamental solution $\wt Y(x,\la)$ (normalized by $\wt Y(0,\la)=E$) of the transformed Hamiltonian system
\eqref{H15} is given by the formula
 \begin{align}& \label{W4}
\wt Y(x,\la)=w_A(x,\la)Y(x,\la)E^*w_A(0,\la)^{-1}E.
\end{align}
Let us set
 \begin{align}& \label{W5}
\clu(\la)=\{\clu_{ij}(\la)\}_{i,j=1}^2:=E^*w_A(0,\la)E,
\end{align}
where $\clu_{ij}(\la)$ are $r\times r$ blocks of $\clu$.
In view of \eqref{W4} and \eqref{W5}, the Weyl circle  (of matrices $\wt M(\la)$)  for the transformed system on $[0,\, \ell']$ and for $\la \in \BC_+$ is determined by
the inequality
 \begin{align}& \nn
\I\begin{bmatrix}I_r & \wt M(\la)^*
\end{bmatrix}\big(\clu(\la)^{-1}\big)^*Y(\ell',\la)^*w_A(x,\la)^*Jw_A(x,\la)Y(\ell',\la)\clu(\la)^{-1} \begin{bmatrix}I_r \\  \wt M(\la)
\end{bmatrix}
\\ & \label{W6}
\leq 0.
\end{align}
Relations \eqref{H10}  and \eqref{H11} yield the following identity \cite[f-la (1.88)]{SaSaR}:
 \begin{align}& \nn
\I w_A(x,\la)^*Jw_A(x,\la)
\\ & \label{W7}
=\I J+\I (\la -\ov \la)\Pi(x)^*(A^*-\ov \la I_n)^{-1}S(x)^{-1} (A- \la I_n)^{-1}\Pi(x).
\end{align}
In this section we consider Hamiltonian systems and assume that $S(0)>0$. Hence, according
to Remark \ref{RkS} we have $S(x)>0$. Now, it is immediate from \eqref{W7} that 
 \begin{align} & \label{W8}
\I w_A(x,\la)^*Jw_A(x,\la)
\leq
\I J \quad (\la \in \BC_+).
\end{align}
Formula \eqref{W8} implies that
 \begin{align}& \nn
\I\begin{bmatrix}I_r & \wt M(\la)^*
\end{bmatrix}\big(\clu(\la)^{-1}\big)^*Y(\ell',\la)^*w_A(x,\la)^*Jw_A(x,\la)Y(\ell',\la)\clu(\la)^{-1} \begin{bmatrix}I_r \\  \wt M(\la)
\end{bmatrix}
\\ & \label{W9}
\leq \I\begin{bmatrix}I_r & \wt M(\la)^*
\end{bmatrix}\big(\clu(\la)^{-1}\big)^*Y(\ell',\la)^*JY(\ell',\la)\clu(\la)^{-1} \begin{bmatrix}I_r \\  \wt M(\la)
\end{bmatrix}
\end{align}
Using \eqref{W9}, we derive the next theorem.
\begin{Tm} \label{TrofWT} Let Hamiltonian system \eqref{E0} $($such that \eqref{E0!} and \eqref{W1} are valid$)$ be given.
Let its GBDT transformation be determined by the triple of matrices $\{A, S(0), \Pi(0)\}$ such that  $S(0)>0$ and that the
matrix identity 
 \begin{align} & \label{W10}
AS(0)-S(0)A^*=\Pi(0)J\Pi(0)^*
\end{align}
holds. Assume that $M(\la)\,\, (\la \in \BC_+)$ belongs to the Weyl
circle \eqref{W2} of the system \eqref{E0} and that
 \begin{align} & \label{W11}
\det \big(\clu_{11}(\la)+\clu_{12}(\la)M(\la)\big)\not= 0,
\end{align}
where $\clu$ is defined in \eqref{W5}. Then 
 \begin{align} & \label{W12}
\wt M(\la)=\big(\clu_{21}(\la)+\clu_{22}(\la)M(\la)\big)\big(\clu_{11}(\la)+\clu_{12}(\la)M(\la)\big)^{-1}
\end{align}
belongs to the Weyl circle of the transformed system.
\end{Tm}
\begin{proof}. Taking into account \eqref{W11} and \eqref{W12}, we obtain
 \begin{align} & \label{W13}
 \begin{bmatrix}I_r \\  \wt M(\la)
\end{bmatrix}
= \clu(\la)\begin{bmatrix}I_r \\   M(\la)
\end{bmatrix}
\big(\clu_{11}(\la)+\clu_{12}(\la)M(\la)\big)^{-1}.
\end{align}
Now, substitute \eqref{W13} into the right-hand side of \eqref{W9} and use \eqref{W2} in order to see that \eqref{W6} is valid.
\end{proof}
According to \cite[p. 671]{Kra}, we have $\I \big(M(\la)-M(\la)^*\big)\leq 0$. Moreover, we have
 \begin{align} & \label{W14}
\I \big(M(\la)-M(\la)^*\big) < 0,
\end{align}
if only $\int_0^{\ell'}y(x,\la)^*H_1(x)y(x,\la)dx>0$ for each nontrivial solution $y$ of \eqref{E0}.
\begin{Rk} \label{dendet} 
If \eqref{W14} is valid, then the inequality \eqref{W11} holds automatically.
Indeed, if $\det \big(\clu_{11}(\la)+\clu_{12}(\la)M(\la)\big)= 0$, then there is a vector $f\not=0$ such that $\big(\clu_{11}(\la)+\clu_{12}(\la)M(\la)\big)f= 0$.
Therefore, recalling that $J$ has the form \eqref{W1}, we obtain
 \begin{align} &\label{W15}
\I f^* \begin{bmatrix}I_r  & M(\la)^*
\end{bmatrix}\clu(\la)^*J
 \clu(\la)\begin{bmatrix}I_r \\   M(\la)
\end{bmatrix}f=0 \quad (f\not=0).
\end{align}
On the other hand, relations \eqref{W5} and \eqref{W8} $($together with the properties of $E$  from \eqref{W3}$)$ imply that
$\clu(\la)^*J\clu(\la)\leq \I J$. Hence, using \eqref{W14}, we derive
 \begin{align} &\label{W16}
\I  \begin{bmatrix}I_r  & M(\la)^*
\end{bmatrix}\clu(\la)^*J
 \clu(\la)\begin{bmatrix}I_r \\   M(\la)
\end{bmatrix}\leq \I  \begin{bmatrix}I_r  & M(\la)^*
\end{bmatrix}J
 \begin{bmatrix}I_r \\   M(\la)
\end{bmatrix}<0,
\end{align}
which contradicts \eqref{W15}.
\end{Rk}
In the limit point case (see, e.g., the discussions in \cite{Kra1, HiSchn})   there is a unique holomorphic in $\BC_+$
Weyl function $\clm(\la)$ the values of which belong to all the Weyl circles \eqref{W2} such that $\ell'<\ell$ ($\la \in \BC_+$).
We note that $\clm(\la)$  is the limit of the values of $M(\la)$ when $\ell'$ tends to $\ell$. Thus,
formula \eqref{W12} shows that
 \begin{align} & \label{W17}
\wt \clm(\la):=\big(\clu_{21}(\la)+\clu_{22}(\la)\clm(\la)\big)\big(\clu_{11}(\la)+\clu_{12}(\la)\clm(\la)\big)^{-1}
\end{align}
is a Weyl function of the transformed system considered on $[0,\ell)$.
%%%%%%%%%%%%%%%%%%%%%%%%%%%%%%%%%%%%%%%%%%%%%%%%%%%%%%%%%
%%%%%%%%%%%%%%%%%%%%%%%%%%%%%%%%%%%%%%%%%%%%%%%%%%%
%%%%%%%%%%%%%%%%%%%%%%%%%%%%%%%%%%%%%%%%
 \section{GBDT  for Shin-Zettl systems}\label{sec2}
 \setcounter{equation}{0}
 Shin-Zettl systems \eqref{E1} present (as well as Hamiltonian systems) an important subclass of systems \eqref{H1}.
 Matrices $Q_1$ and $Q_2$, in the case of  Shin-Zettl systems, have the form
  \begin{align}& \label{E8}
Q_1(x)= \begin{bmatrix} 0 & 0 \\ \om(x) & 0
\end{bmatrix}, \quad Q_0(x)=-
\begin{bmatrix} r_1(x) & p(x)^{-1} \\ q(x)  & \, r_2(x)
\end{bmatrix}.
\end{align}
 Recall that GBDT is determined by the parameter matrices $A_1$, $A_2$, $S(0)$, $\Pi_1(0)$ and  $\Pi_2(0)$ such
 that \eqref{E6} holds.  For  the Shin-Zettl systems, we have $m=2$, and so matrices $\Pi_1(0)$ and  $\Pi_2(0)$
 are $n \times 2 $ matrices. Using   the second equality in \eqref{E1} and the first equality in \eqref{E8}, we rewrite $\wt F$ given
by \eqref{E11}--\eqref{E12} in the Shin-Zettl form
\begin{align}& \label{E14}
\wt F(x,\la)=\begin{bmatrix}\wt r_1(x) & \wt p(x)^{-1} \\ \wt q(x) - \la \wt \om(x) & \, \wt r_2(x)
\end{bmatrix}, \quad \wt \om=\om, \quad \wt p=p;
\\ & \label{E15}
\wt r_1=r_1-\om X_{12}, \quad \wt r_2=r_2+\om X_{12}, \quad \wt q=q+\om(X_{11}-X_{22}).
\end{align}
where $X_{ik}(x)$ are the entries of $X(x)$. Now, the following proposition is immediate from
Theorem \ref{Tm1GBDT}.
\begin{Pn}\label{PnGBDT} Let $y(x,\la)$ satisfy Shin-Zettl system \eqref{E1} and let
$w_A$ be given by \eqref{E10}, where the matrix functions $\Pi_1$, $\Pi_2$ and $S$ are determined
by \eqref{E7} and identity \eqref{E6} holds. Then the function $\wt y(x,\la)=w_A(x,\la)y(x,\la)$
satisfies, in the points of invertibility of $S(x)$,  the transformed Shin-Zettl system \eqref{H4},
where $\wt F(x,\la)$ is given by \eqref{E14} and \eqref{E15}.
\end{Pn}

%%%%%%%%%%%%%%%%%%%%%%%%%%%
The next corollary easily follows from Proposition \ref{PnGBDT}.
\begin{Cy} \label{CyJSym} Let the conditions of Proposition \ref{PnGBDT} hold
and let the initial system \eqref{E1} be Lagrange-J-symmetric $($i.e., let \eqref{E5} be valid$)$. Then the transformed system
is Lagrange-J-symmetric as well, that is, the equality $ \wt r_1=-\wt r_2$ holds.
\end{Cy}
In the next section, we consider the Lagrange-symmetric case (i.e., the case \eqref{E4}).
%%%%%%%%%%%%%%%%%%%%%%%%%%%%%%%%%%%
%%%%%%%%%%%%%%%%%%%%%%%%%%%%%%%%%%
 \section{Lagrange-symmetric case}\label{sec3}
 \setcounter{equation}{0}
 Further we assume that \eqref{E4} is fulfilled and rewrite \eqref{E8} for that case:
 \begin{align}& \label{E8!}
Q_1(x)= \begin{bmatrix} 0 & 0 \\ \om(x) & 0
\end{bmatrix}, \quad Q_0(x)=-
\begin{bmatrix}
 r(x) & p(x)^{-1} \\ q(x)  & \, -\ov{r(x)}
\end{bmatrix}, \\
& \label{E8!!}
 r(x):=r_1(x)=-\ov{r_2(x)}.
\end{align}
Now, system
 \eqref{E1} may be rewritten as the quasi-differential equation:
\begin{align}& \label{L0}
-\big(u^{[1]}\big)^{\prime}-\ov{r}u^{[1]}+qu=\la \om u, \qquad u^{[1]}:=p(u^{\prime}-ru),
\end{align} 
where $y_1(x)=u(x)$, $y_2(x)=u^{[1]}(x)$ and the quasi-differential expression
$$Mu=-\big(u^{[1]}\big)^{\prime}-\ov{r}u^{[1]}+qu-\la \om u$$
is symmetric (see, e.g., \cite{Ev2}). See also \cite{BPT, Nai, Z} and references therein
on symmetric expressions  $(-\big(u^{[1]}\big)^{\prime}-\ov{r}u^{[1]}+qu)/\om$
in the weighted spaces $L^2_{|\om|}[0,\ell)$ and $L^2_{\om}[0,\ell)$.
Using the quasi-derivative $u^{[1]}$ one may consider
Sturm-Liouville equations  (including self-adjoint Sturm-Liouville equations) 
with non-smooth coefficients (see, e.g., the discussions in \cite[p. 455]{Ze} and in
\cite[p. 25]{Z}).

We note that $Q_1$ and $Q_0$ given by \eqref{E8!} admit representation \eqref{H5},
where 
 \begin{align}& \label{H17}
J=\I \s_2, \quad H_1(x)= \begin{bmatrix}  \om(x) & 0\\ 0 &0
\end{bmatrix}, \quad H_0(x)=
\begin{bmatrix}
  - q(x)  & \, \ov{r(x)} \\ r(x) & p(x)^{-1}
\end{bmatrix},
  \end{align}
$\s_2:=\begin{bmatrix}0 & -\I \\ \I & 0
\end{bmatrix}$ is a  Pauli matrix, and  \eqref{H6} holds. In fact, conditions \eqref{H5} and \eqref{H6}
(in the Shin-Zettl  case and with $J=\I \s_2$)
are equivalent to the conditions \eqref{E4} of Lagrange symmetry. (Clearly, when $\om \geq 0$ we deal
with a subclass of Hamiltonian systems.) 
 Thus, omitting the indices in $A_1$ and $\Pi_1$ and rewriting \eqref{H7} in the form
 \begin{align}& \label{L1}
A=A_1, \quad \Pi=\Pi_1; \quad  A_2=A^*, \quad S(0)=S(0)^*, \quad \Pi_2(0)=-\I  \Pi(0)\s_2,
 \end{align}
 we see that the formulas of $\S$2 in Section \ref{sec1'} are valid for Lagrange-symmetric case.
 
Since $\Pi_2(x)=-\I \Pi(x)\s_2$, formula
 \eqref{E12} for $X$ may be rewritten as
\begin{align}& \label{L6}
X(x)=J\Pi(x)^*S(x)^{-1}\Pi(x), \quad J=\I \s_2=\begin{bmatrix}0 & 1 \\ -1 & 0
\end{bmatrix},
\end{align} 
and we obtain
\begin{align}& \label{L7}
X_{12}(x)=\ov{X_{12}(x)}, \quad X_{22}(x)=-\ov{X_{11}(x)}.
\end{align} 
Recall that $\Pi(x)$ and $S(x)$ are given by the equations
 \begin{align}& \label{L3'}
\Pi^{\prime}=-A\Pi J H_1- \Pi J H_0, \quad  S^{\prime}=\Pi J H_1J^* \Pi^*, \quad J=\I \s_2.
\end{align}
Formula \eqref{H11} for the Darboux matrix takes the form
 \begin{align}
& \label{L5}
w_A(x,\la)=I_2-\I \s_2 \Pi(x)^*S(x)^{-1}(A-\la I_n)^{-1}\Pi(x).
\end{align}
Using \eqref{H17} and \eqref{L7}, we derive from Propositions \ref{Cy1GBDT} and \ref{PnGBDT} the next corollary.
\begin{Cy} \label{CySym} Assume that the initial Shin-Zettl system is Lagrange-symmet-ric $($i.e., that \eqref{E4} holds$)$.
Let the matrices $A$, $\Pi(0)$ and $S(0)$ be chosen so that $S(0)=S(0)^*$ and $AS(0)-S(0)A^*=\I \Pi(0)\s_2 \Pi(0)^*$,
and let $\Pi(x)$, $S(x)$ and $X(x)$ be determined by \eqref{L3'} and \eqref{L6}, respectively.

Then the corresponding transformed Shin-Zettl system \\ $\wt y^{\prime}(x,\la)= \wt F(x,\la)
\wt y(x,\la)$
is given by \eqref{E14}, where
 \begin{align}
& \label{Z1}
\wt r_1=-\ov{\wt r_2}=r-\om X_{12}, \quad \wt q=q+\om (X_{11}+\ov{X_{11}}).
\end{align}
This transformed system is 
Lagrange-symmetric as well. Moreover, the function 
$\wt y(x,\la)=w_A(x,\la)y(x,\la)$, 
where $w_A$ has the form
\eqref{L5}, satisfies the transformed system.
\end{Cy}

\section{Sturm-Liouville equations}\label{StL}
\setcounter{equation}{0}
In this section we consider Sturm-Liouville equation \eqref{E3}. GBDT for its particular case (namely, for Schr\"odinger equation
where $p\equiv \om \equiv 1$) was dealt with in \cite{GKS3} but the general equation \eqref{E3} contains other interesting subcases,
where GBDT could be useful as well.
\begin{Pn} \label{Pn4.1} Let the function $p\om$ be differentiable and its derivative $(p\om)^{\prime}$ as well as the functions $p^{-1}$, $q$
and $\om$ be locally summable on $[0, \, \ell)$. Assume that
\begin{align}& \label{E4!}
\om= \ov \om, \quad p=\ov p, \quad q=\ov q, \quad r \equiv 0,
\end{align}
and set
\begin{align}& \label{St1}
\wt y(x,\la)=w_A(x,\la)y(x,\la),
\end{align}
where $w_A$ is given by the relations \eqref{L5} and \eqref{L3'}, $H_0$ and $H_1$ $($in \eqref{L3'}$)$ are given by \eqref{H17}
and $y$ satisfies the initial Lagrange-symmetric Shin-Zettl equation
\begin{align}& \label{E1!}
y^{\prime}(x,\la)= J\big(\la H_1(x)+ H_0(x)\big)y(x, \la).
\end{align}
Then the entry $\wt y_1$ of $\wt y$  satisfies the transformed Sturm-Liouville equation
\begin{align}& \label{E3!}
-\big( p(x) \wt y_1^{\prime}(x,\la)\big)^{\prime}+ \breve q(x)\wt y_1(x,\la)=\la \om(x) \wt y_1(x,\la),
\end{align}
where
\begin{align}& \label{St2}
\breve q=q+2\om (X_{11}-X_{22})+2p(\om X_{12})^2-(p\om)^{\prime}X_{12},
\end{align}
and $X_{ik}$ are the entries of $X$ given by \eqref{L6}.
\end{Pn}
\begin{proof}. Recall that in Section \ref{sec3} we rewrote Shin-Zettl system in the form \eqref{L0} where $u=y_1$.
In the notations of the transformed system it means
\begin{align}& \label{St3}
-\big(p(\wt y_1^{\prime}-\wt r \wt y_1)\big)^{\prime}-\ov{\wt r}p(\wt y_1^{\prime}-\wt r \wt y_1)+\wt q \wt y_1=\la \om \wt y_1,
\end{align} 
where $\wt r:=\wt r_1(x)=-\ov{\wt r_2(x)}$. Using the  identity $r_1\equiv r_2\equiv 0$ and equalities \eqref{E15} and \eqref{L7} we present \eqref{St3} in the form
\begin{align}& \nonumber
-\big(p\wt y_1^{\prime}\big)^{\prime}-\big(p\om X_{12}\wt y_1\big)^{\prime}+p\om X_{12}\wt y_1^{\prime}+p(\om X_{12})^2\wt y_1
\\ & \label{St4}
+\big(q+\om(X_{11}-X_{22})\big) \wt y_1=\la \om \wt y_1,
\end{align} 
which is equivalent to \eqref{E3!} with
\begin{align}& \label{St5}
\breve q=q+\om (X_{11}-X_{22})+p(\om X_{12})^2-(p\om)^{\prime}X_{12}-p\om X_{12}^{\prime}.
\end{align}
Finally, in order to show that  the functions $\breve q$ given by \eqref{St2} and \eqref{St5} coincide,
let us differentiate $X_{12}$. Taking into account \eqref{L6} and \eqref{L3'},  we obtain:
\begin{align}& \nn
X^{\prime}=J (H_1J\Pi^*A^*S^{-1}\Pi+H_0 J\Pi^*S^{-1}\Pi)-J\Pi^*S^{-1}\Pi J H_1J^* \Pi^*S^{-1}\Pi
\\ & \nn
-J\Pi^*S^{-1}(A\Pi J H_1+\Pi J H_0).
\end{align}
In particular, for $X_{12}$ we obtain 
\begin{align}& \label{St6}
X_{12}^{\prime}= p^{-1}(X_{22} -X_{11})-\om X_{12}^2.
\end{align}
Here we again took into account that $r_1\equiv r_2\equiv 0$. Equalities  \eqref{St5} and \eqref{St6} imply \eqref{St2}.
\end{proof}
\begin{Rk} \label{Rkrv} In view of \eqref{L7}, \eqref{E4!} and \eqref{St2}, the equality $\Im(\breve q)\equiv 0$ is valid.
Thus, the coefficients of the transformed Sturm-Liouville equation \eqref{E3!} are real-valued. It is easy to see
that  the function $\breve q$ is locally summable on $[0,\, \ell)$ if the conditions of Proposition \ref{Pn4.1} hold
and $S(x)$ is invertible on $[0,\, \ell)$.
\end{Rk}
%%%%%%%%%%%%%%%%%%%%%%%%%%%%%%%%%%%%%%%%%%%%%%
%%%%%%%%%%%%%%%%%%%%%%%%%%%%%%%%%%%%%%%%%%%%%%%%
\section{Dynamical systems}\label{DS}
\setcounter{equation}{0}
%%%%%%%%%%%%%%%%%%%%
\subsection{Dynamical symplectic system}\label{DSympl}
Formally applying Laplace transform  to the  system \eqref{E0} (satisfying \eqref{H6}), we come to the
interesting dynamical system
\begin{align}& \label{D1}
\frac{\p}{\p x}z(x,t)=J\left(- H_1(x)\frac{\p}{\p t}z(x,t)+H_0(x)z(x,t)\right).
\end{align}
When $J^*=J^{-1}$ system \eqref{D1} is a dynamical symplectic system.

In order to construct Darboux transformation of system \eqref{D1} and solutions of the transformed system, we use \eqref{H9'} and rewrite \eqref{D0} (for our
case where the relations \eqref{H5}--\eqref{H7} are valid) in the form
\begin{align}& \label{D2}
\big(J\Pi^*S^{-1}\big)^{\prime}=J\big(H_1J \Pi^*S^{-1}A+\wt H_0 J\Pi^*S^{-1}\big),\\
& \label{D3}
\wt H_0=H_0-X^*H_1-H_1 X, \quad X=J \Pi^*S^{-1}\Pi .
\end{align}
We note that \eqref{D3} is equivalent to the second equality in \eqref{H12}.
%%%%%%%%%%%%%%%%%%%%%%%%%%%%%%%%%%
\begin{Pn}\label{PnDSympl} Let $J$, $H_1(x)$ and $H_0(x)$ satisfying    \eqref{H6}, as well as 
the triple $\{A, S(0)=S(0)^*, \Pi(0)\}$ satisfying \eqref{W10}, be given.  
Let
the matrix functions $\Pi(x)$  and $S(x)$ be determined by \eqref{H9}. Then the vector functions
\begin{align}& \label{D4}
\wt z(x,t)=J\Pi(x)^*S(x)^{-1} \E^{-t A}h  \quad (h\in \BC^m)
\end{align}
satisfy, in the points of invertibility of $S(x)$, the transformed  dynamical   system $($of the same form as \eqref{D1}$)$. 
More precisely, we have
\begin{align}& \label{D5}
\frac{\p}{\p x}\wt z(x,t)=J\left(- H_1(x)\frac{\p}{\p t}\wt z(x,t)+\wt H_0(x)\wt z(x,t)\right),
\end{align}
where $\wt H_0$ is given by \eqref{D3}
\end{Pn}
\begin{proof}.
In view of \eqref{D2} and \eqref{D4}, both sides of \eqref{D5} equal
$J\big(H_1J\Pi^*S^{-1}A+\wt H_0 J \Pi^*S^{-1}\big)\E^{-tA}h$.
\end{proof}
When $H_1 \geq 0$, the energy $E_{z}(t)$ of the solutions $z$ of  system \eqref{D1} on $[0, a]$ ($0<a<\ell$)
is given by the formula
\begin{align}& \label{D6-}
E_{z}(t)^2=\int_0^a z(x,t)^*H_1(x)z(x,t)dx.
\end{align}
The energy of the transformed solutions $\wt z$ of the form \eqref{D4} is expressed via $A$ and $S(x)$.
\begin{Pn}\label{PnE} Let the conditions of Proposition \ref{PnDSympl} hold and assume additionally that $H_1 \geq 0$ and $S(0) > 0$.
Then the energy $E_{\wt z}$, where $\wt z$ has the form \eqref{D4}, is  given by the formula
\begin{align}& \label{D7-}
E_{\wt z}(t)=\sqrt{h^*\E^{- tA^*}\big(S(0)^{-1}-S(a)^{-1}\big)\E^{- tA}h}.
\end{align}
\end{Pn}
\begin{proof}. Taking into account \eqref{L3'} and \eqref{D4} we see that
\begin{align}& \label{D8-}
 \wt z(x,t)^*H_1(x)\wt z(x,t)=-h^*\E^{- tA^*}\big(S(x)^{-1}\big)^{\prime}\E^{- tA}h.
\end{align}
Formula \eqref{D7-} follows from \eqref{D6-} and \eqref{D8-}.
\end{proof}

%%%%%%%%%%%%%%%%%%%%%%%%%%%%%%%%%%%%%%%%%%%%%%%%%%%%%%%%%%%%%
\subsection{Two-way diffusion equation}\label{DifEq}
In this subsection, we consider the important case when 
$J$, $H_1$ and $H_0$  have the form \eqref{H17} (i.e., the same form as in Lagrange-symmetric Shin-Zettl system)
and $\om= \ov \om$,  $p=\ov p$,   $q=\ov q$.
In that case we set 
\begin{align}& \label{D6}
z(x,t)=\begin{bmatrix}z_1(x,t) \\ z_2(x,t) \end{bmatrix}, \quad \wt z(x,t)=\begin{bmatrix}\wt z_1(x,t) \\  \wt z_2(x,t) \end{bmatrix},
\end{align}
and rewrite \eqref{D1} in the form
\begin{align}& \label{D7}
z_1^{\prime}=r z_1+p^{-1}z_2, \quad z_2^{\prime}=\om \frac{\p}{\p t}z_1+qz_1-\ov{r}z_2 \quad \Big(z_i^{\prime}=\frac{\p}{\p x}z_i\Big).
\end{align}
Next, we rewrite the first equality in \eqref{D7} as $z_2=p(z_1^{\prime}-r z_1)$, substitute the expression for $z_2$ into
the second equality in \eqref{D7} and obtain
\begin{align}& \label{D8}
\om \frac{\p}{\p t}z_1=\big(p(z_1^{\prime}-r z_1)\big)^{\prime} -qz_1+\ov{r}p(z_1^{\prime}-r z_1).
\end{align}
In particular, when $r=0$, equation \eqref{D8} takes the form
 \begin{align}& \label{D10}
\om \frac{\p}{\p t}z_1=\big(p(z_1^{\prime})\big)^{\prime} -qz_1.
\end{align}
We note that equation \eqref{D10} coincides (in the case of sign-indefinite $\om$) with the 
{\it two-way diffusion equation} (6.1) in \cite{Kost}. See also various references in \cite{Bea, FiK, Kost}
on the literature related to the two-way diffusion equation.

According to Corollary \ref{CySym}, $\wt H_0$ has the same form as $H_0$.
More precisely, we have (see \eqref{Z1} or \eqref{D3}):
\begin{equation} \label{D9}
\wt H_0(x)=\begin{bmatrix} -\wt q(x) & \ov{\wt r(x)}
\\ \wt r(x) & p(x)^{-1}
\end{bmatrix}, \quad \wt r=r-\om X_{12}, \quad \wt q=q+\om(X_{11}+\ov{X_{11}}).
\end{equation}
In the same way as \eqref{D1} yields \eqref{D8}, equation \eqref{D5} implies that
 the entry $\wt z_1$ of the solution $\wt z$ given by  \eqref{D4} satisfies the equation
\begin{align}& \label{D11}
\om \frac{\p}{\p t}\wt z_1=\big(p(\wt z_1^{\prime}-\wt r \wt z_1)\big)^{\prime} -\wt q \wt z_1+\ov{\wt r}p(\wt z_1^{\prime}-\wt r \wt z_1).
\end{align}
Assuming $r\equiv 0$, we see that
\begin{align} 
& \label{D12}
\wt r=-\om X_{12}, \quad \wt q=q+\om(X_{11}+\ov{X_{11}}).
\end{align}
Multiplying the left-hand side of \eqref{St3} by ``$-1$" and substituting there $\wt y_1=\wt z_1$ we obtain
the right-hand side of \eqref{D11}. Hence, the proof of Proposition of \ref{Pn4.1} shows that
\begin{align}& \label{D13}
\big(p(\wt z_1^{\prime}-\wt r \wt z_1)\big)^{\prime} -\wt q \wt z_1+\ov{\wt r}p(\wt z_1^{\prime}-\wt r \wt z_1)
= \big( p\wt z_1^{\prime}\big)^{\prime}- \breve q\wt z_1,
\\ & \label{D14}
\breve q=q+2\om (X_{11}-X_{22})+2p(\om X_{12})^2-(p\om)^{\prime}X_{12}.
\end{align}
From \eqref{D6}, \eqref{D11} and \eqref{D13}, the next proposition is immediate.
\begin{Pn}\label{PnTwoW} Let $J$, $H_1$ and $H_0$  have the form \eqref{H17},
and let the function $p\om$ be differentiable and its derivative $(p\om)^{\prime}$ as well as the functions $p^{-1}$, $q$
and $\om$ be locally summable on $[0, \, \ell)$. Assume that \eqref{E4!} holds and that the triple $\{A, S(0)=S(0)^*, \Pi(0)\}$ satisfies \eqref{W10}.
Introduce $\Pi(x)$ and $S(x)$ via
 \eqref{H9}.
 
 Then the function $\wt z_1$  $($given by
\eqref{D4} and \eqref{H9}$)$ satisfies, in the points of invertibility of $S(x)$, the dynamical equation
\begin{align}& \label{D15}
\om \frac{\p}{\p t}\wt z_1=\big( p\wt z_1^{\prime}\big)^{\prime}- \breve q\wt z_1,
\end{align}
where $\breve q$  is given by \eqref{D14}.
\end{Pn}
Recall that \eqref{D15} is an equation of the form \eqref{D10}.
\begin{Rk}\label{RkInt} It is important that \cite[Theorem 7.4]{SaSaR} and our Theorem \ref{Tm1GBDT}, in particular, is valid
on any interval $\cli$ such that $0 \in \cli$. Thus, the previous statements of the paper, excluding the last sentence
in Proposition \ref{Cy1GBDT}, Remark \ref{RkS}, Proposition \ref{PnE} and
the statements from
Section \ref{secAp}   $($where the condition $S(x)>0$ is essential$)$, 
are also valid on the intervals $\cli$ such that $0 \in \cli$. The interval $[0, \ell)$ was chosen for simplicity but
the interval $(-\ell, \ell)$ is sometimes more convenient in the  two-way diffusion equation and in the indefinite Sturm-Liouville case. 
\end{Rk}
%%%%%%%%%%%%%%%%%%%%%%%%%%%%%%%%%%%%%%%%%%%%%%
%%%%%%%%%%%%%%%%%%%%%%%%%%%%%%%%%%%%%%%%%%%%%%%%%
\section{Indefinite Sturm-Liouville equations}\label{IStL}
\setcounter{equation}{0}
Symplectic systems and indefinite Sturm-Liouville equations are of growing interest in the literature
(see, e.g., \cite{AtkEv, BiL, BPT, LaW, Levi, Qi} and references therein). Therefore, in this  section
of the paper we shall consider some examples of Darboux transformation for the Lagrange-symmetric Shin-Zettl system and
indefinite Sturm-Liouville equation
considered on the interval $(-\ell, \ell)$, see Remark \ref{RkInt}.

More precisely, we shall construct explicit solutions for the interesting model case
\begin{align}& \label{I1}
\om(x)={\mathrm{sgn}}(x), \quad p(x) \equiv 1,
\end{align}
which was studied in \cite{KarT}. First, we consider Shin-Zettl system \eqref{E0}, \eqref{H17} and assume
that the equalities \eqref{I1} and 
\begin{align}& \label{I2}
q(x)=r(x)\equiv 0
\end{align}
hold for the initial system.  We consider Darboux transformations
determined by the triples of matrices $\{A,\, S(0), \, \Pi(0)\}$ of the form
\begin{align}& \label{I3}
A=\a^2, \quad S(0)=0, \quad \Pi(0)=\begin{bmatrix} -2\I \a g & 2\mu \a g \end{bmatrix},
\end{align}
where $\a$ are $n \times n$ matrices, $g \in \BC^n$ are vector columns, $\mu$
are purely imaginary values (i.e. $\ov{\mu}=-\mu$), and
\begin{align}& \label{I4}
\det(\mu \a \pm I_n)\not=0, \quad \det(\mu \a \pm \I I_n)\not=0.
\end{align}
It is easily checked that the third equality in \eqref{I3} yields
$\Pi(0)J\Pi(0)^*=0$ $\,(J=\I \s_2)$, and so the matrix identity \eqref{W10},
which is required in GBDT,  holds for the triple of the
form \eqref{I3}.

Next, we partition $\Pi(x)$ into two columns $\Pi(x)=\begin{bmatrix}\Lam_1(x) &\Lam_2(x) \end{bmatrix}$,
and (taking into account \eqref{I1}--\eqref{I3}) rewrite the first system in \eqref{L3'} in the form
\begin{equation} \label{I5}
\Lam_1^{\prime}= \left\{ \begin{array}{lr} \a^2 \Lam_2 &  {\mathrm{for}} \quad x>0 \\
-\a^2 \Lam_2 & {\mathrm{for}} \quad x<0 \end{array}  \right. ; \qquad \Lam_2^{\prime}=-\Lam_1.
\end{equation}
It is immediate that the vector functions
\begin{align}& \label{I6}
\begin{array}{lr} \Lam_1(x)=-\I \a \big(\E^{\I x \a}(\mu \a +I_n)g-\E^{-\I x \a}(\mu \a -I_n)g\big) \\
\Lam_2(x)= \E^{\I x \a}(\mu \a +I_n)g+ \E^{-\I x \a}(\mu \a -I_n)g \end{array} \quad {\mathrm{for}} \quad x \geq 0;
\\  & \label{I7}
\begin{array}{lr} \Lam_1(x)=- \a \big(\E^{x \a}(\mu \a +\I I_n)g-\E^{- x \a}(\mu \a - \I I_n)g\big) \\
\Lam_2(x)= \E^{ x \a}(\mu \a +\I I_n)g+ \E^{- x \a}(\mu \a -\I I_n)g \end{array} \quad {\mathrm{for}} \quad x \leq 0
\end{align}
satisfy \eqref{I5} and the third equality in \eqref{I3}.

The second system in  \eqref{L3'} takes the form $S^{\prime}=\om \Lam_2 \Lam_2^*$.
Hence, using $S(0)=0$, we see that
\begin{align}& \label{I8}
S(x)=\int_0^x \Lam_2(t) \Lam_2(t)^*dt \geq 0\quad  (x >0), \\
& \label{I9}
 S(x)=\int_x^0 \Lam_2(t) \Lam_2(t)^*dt \geq 0 \quad (x <0).
\end{align}
\begin{Rk} \label{S>0} It follows from \eqref{I8} and \eqref{I9} that usually we have 
\begin{align}& \label{I10}
S(x)>0  \quad  {\mathrm{for}} \quad x\not=0.
\end{align}
In particular, \eqref{I10} holds when the pair $\{\wh \a, \, \wh g\}$, where 
$$\wh \a:= \begin{bmatrix} \a & 0
\\ 0 & - \a
\end{bmatrix}, \quad \wh g:= \begin{bmatrix} g
\\ g
\end{bmatrix},$$
is controllable. Indeed, if \eqref{I10} is not valid, then there is $f\in \BC^n$ $\,(f\not=0)$ such that
$f^*\Lam_2(x)=0$ either for all $x>0$ or for all $x<0$. In view of \eqref{I4}, \eqref{I6} and \eqref{I7}
it means that $\wh f^*\E^{cx\wh \a}\wh g\equiv 0$ for some $\wh f \in \BC^{2n}$ and $c\in \BC$ $\,(\wh f \not=0, \quad c\not=0)$.
However, this contradicts the controllability of $\{\wh \a, \, \wh g\}$ $($see, e.g., \cite{Cop}$)$.
\end{Rk}
Formulas \eqref{I6}--\eqref{I9} present explicit expressions for $\Pi(x)$ and $S(x)$, and so the Darboux matrix $w_A(x,\la)$
of the form \eqref{L5} is constructed explicitly. 

In order to use Corollary \ref{CySym} we also solve explicitly the initial
Shin-Zettl system \eqref{E0}, \eqref{H17}, where \eqref{I1} and \eqref{I2} hold. Namely, we introduce matrices
\begin{align}& \label{I11}
T_+(\la)=\begin{bmatrix}1 & 1 \\ 
\I \sqrt{\la} & -\I \sqrt{\la}
 \end{bmatrix}, \quad
D_+(\la)=\begin{bmatrix}
\I \sqrt{\la} & 0 \\ 0 & -\I \sqrt{\la}
 \end{bmatrix} ,
\\ & \label{I12} 
T_-(\la)=\begin{bmatrix}1 & 1 \\ 
\sqrt{\la} & - \sqrt{\la}
 \end{bmatrix}, \quad
D_-(\la)=\begin{bmatrix}
\sqrt{\la} & 0 \\ 0 & - \sqrt{\la}
 \end{bmatrix} ,
\end{align}
where $\sqrt{\la}$ is any fixed branch of the square root of $\la$. It is easy to see that (in our case) $F$ given
in \eqref{E0} satisfies the equalities $FT_+=T_+D_+$ for $x>0$ and $FT_-=T_-D_-$ for $x<0$.
Therefore, solutions $y$ of the initial Shin-Zettl system \eqref{E0} are given by the formulas
\begin{align}& \label{I13}
y(x,\la)=T_+(\la)\E^{x D_+(\la)} T_+(\la)^{-1}h \quad (x>0), 
\\ & \label{I13'}
 y(x,\la)=T_-(\la)\E^{x D_-(\la)} T_-(\la)^{-1}h \quad (x<0)
\end{align}
with any vectors $h\in \BC^2$. Now, Corollary \ref{CySym} and Remarks \ref{RkInt} and \ref{S>0}
yield our next corollary.
\begin{Cy} \label{CyInd1}
Assume that the initial Shin-Zettl system on  $(-\ell, \ell)$ has the form \eqref{E0}, \eqref{H17}
  and that equalities \eqref{I1} and \eqref{I2}
 hold.
Let the matrices $A$, $\Pi(0)$ and $S(0)$ have the form \eqref{I3}. Then the corresponding GBDT-transformed Shin-Zettl system 
$\wt y^{\prime}(x,\la)= J(\la H_1(x)+\wt H_0(x))
\wt y(x,\la)$ is Lagrange symmetric
and we have $\wt H_0(x)=\begin{bmatrix}
  - \wt q(x)  & \, \ov{\wt r(x)} \\ \wt r(x) & 1
\end{bmatrix}$, where
\begin{align}& \label{I14}
\wt r(x)=-{\mathrm{sgn}}(x)X_{12}(x), \quad \wt q(x)={\mathrm{sgn}}(x)(X_{11}(x)+\ov{X_{11}(x)}),
\end{align}
$X_{ij}$ are the blocks of $X=J\Pi^*S^{-1}\Pi$, and explicit expressions for $S$ and $\Pi$ are
given in \eqref{I6}--\eqref{I9}. The controllability of the pair $\{\wh \a, \, \wh g\}$ is a sufficient condition of the invertibility
of $S(x)$ at $x\not=0$. If, indeed, $\det S(x)\not=0$ for $x\not=0$, then  $X(x)$ and the Darboux matrix $w_A(x,\la)$ of the form \eqref{L5}
are well defined and explicitly expressed via $\Pi(x)$ and $S(x)$ at $x\not=0$. Moreover, the solutions $\wt y$
of the GBDT-transformed system are explicitly expressed via the formula $\wt y(x,\la)=w_A(x,\la)y(x,\la),$ where $y$
is given by \eqref{I13}and \eqref{I13'}.
\end{Cy}
By virtue of Proposition \ref{Pn4.1}, Remark \ref{RkInt} and Corollary \ref{CyInd1}, we obtain explicit
solutions of indefinite Sturm-Liouville systems
\begin{align}& \label{I15}
- \wt y_1^{\prime \prime}(x,\la)+\breve q(x)\wt y_1(x,\la)=\la {\mathrm{sgn}}(x)\wt y_1(x,\la) \quad (-\ell<x<\ell).
\end{align}
\begin{Cy} \label{CyInd2} Let $\Pi(x)$ and $S(x)$ be given by \eqref{I6}--\eqref{I9} and assume that $\det S(x)\not=0$ for $x \not=0$.
Set $\wt y(x,\la)=w_A(x,\la)y(x,\la)$ where explicit expressions for $w_A(x,\la)$ $($with $A=\a^2)$ and $y(x,\la)$ are given by \eqref{L5} and \eqref{I13},
\eqref{I13'},
respectively. Then the first entry $\wt y_1$ of $y$ satisfies the indefinite Sturm-Liouville system \eqref{I15} where
\begin{align}& \label{I16}
\breve q(x)=2{\mathrm{sgn}}(x)\big(X_{11}(x)-X_{22}(x)\big)+2X_{12}(x)^2
\end{align}
and $X_{ij}$ are the blocks of $X=J\Pi^*S^{-1}\Pi$.
\end{Cy}
The singularity of $\breve q(x)$ at $x=0$ is of interest.  Some particular cases (but in greater detail) were considered in  \cite[Section 5]{KoSaTe},
 and it was proved for those cases that $\breve q(x)=O(x^{-2})$ when $x$ tends to $0$.

\bigskip

\noindent{\bf Acknowledgments.}
 {This research   was supported by the
Austrian Science Fund (FWF) under Grant  No. P29177.}
%%%%%%%%%%%%%%%%%%%%%%%%%%%%%%%%%%%%%%%%%%%%%%
%%%%%%%%%%%%%%%%%%%%%%%%%%%%%%%%%%%%%%%%%%%%%%%
\newpage
%%%%%%%%%%%%%%%%%%%%%%%%%%%%%%%%%%%%%%%%%%%
%%%%%%%%%%%%%%%%%%%%%%%%%%%%%%%%%%%%%%%%%%%%

\begin{flushright}

A.L. Sakhnovich,\\
Fakult\"at f\"ur Mathematik, Universit\"at Wien, \\
Oskar-Morgenstern-Platz 1, A-1090 Vienna, Austria
\end{flushright}

%%%%%%%%%%%%%%%%%%%%%%%%

\end{document}